\documentclass[CRMATH,Unicode,manuscript]{cedram}

\usepackage{amsmath,amsfonts,graphicx,tikz,pgfplots, tikz-3dplot, pdfsync}
\usepackage{bm} % boldface greeks

\usepackage{tikz,pgfplots,bm}
\usetikzlibrary{patterns}

\def\eq#1{{\begin{eqnarray}#1\end{eqnarray}}}
\def\eqn#1{{\begin{eqnarray*}#1\end{eqnarray*}}}

\def\R{\mathbb R}
\def \E{{\mathbb{E}}}
\def \N{{\mathbb{N}}}
\def\vn{\mathbf n}
\def\va{\mathbf a}
\def\vq{\mathbf q}
\def\vx{\mathbf x}
\def\vy{\mathbf y}
\def\ds{\displaystyle }
\def\add#1{\color{red}#1\color{black}}

\title{The Dual Characteristic-Galerkin Method}

\alttitle{La méthode des caract\'eristiques-Galerkin duale}

\author{\firstname{Frédéric}  \lastname{Hecht}}
\address{LJLL, Boite 187, Sorbonne Université, Place Jussieu, 75005 Paris, France}
\email[]{frederic.hecht@academie-sciences.fr}
\author{\firstname{Olivier} \lastname{Pironneau}}
\address{LJLL, Boite 187, Sorbonne Université, Place Jussieu, 75005 Paris, France}
\email[]{olivier.pironneau@academie-sciences.fr}


\begin{abstract}
The Dual Characteristic-Galerkin method (DCGM) is conservative, precise and experimentally positive. We present the method and prove convergence and $L^2$-stability in the case of Neumann boundary conditions. In a 2D numerical finite element setting (FEM), the method is compared to  Primal Characteristic-Galerkin (PCGM), Streamline upwinding (SUPG), the Dual Discontinuous Galerkin method (DDG) and centered FEM without upwinding. DCGM is difficult to implement numerically but, in the numerical context of this note, it is far superior to all others.
 \end{abstract}

\begin{altabstract} 
La méthode \emph{Dual Characteristic-Galerkin} (DCGM) est conservative, précise et expérimentalement positive. Nous prouvons la convergence et la stabilité $L^2$. Dans le cadre numérique des méthodes d'éléments finis (FEM) en 2D, la méthode est comparée à la méthode \emph{Primal Characteristic-Galerkin} (PCGM), au \emph{Streamline upwinding} (SUPG), à la méthode \emph{Dual Discontinuous Galerkin} (DDG) et à une discretisation FEM sans décentrage. La méthode DCGM est difficile à mettre en œuvre numériquement, mais elle est de loin supérieure à toutes les autres dans le cadre étudié dans cette note.
\end{altabstract}

\begin{document}
\maketitle

\section*{Introduction}
A good numerical method for the convection-diffusion equation is important in itself but it is also a test bed for more complex systems such as the Navier-Stokes equations. A finite element method (FEM) combined with a first or second order implicit in time discretization without upwinding works only if a CFL condition is satisfied, a severe constraint if the viscous coefficient is small (the method is also known as Arakawa's scheme in meteorology \cite{arak}). Hence in the eighties a number of upwinding schemes have been proposed in particular by K. Baba et al \cite{baba}, J.-P. Benque et al \cite{benque} T.J.R. Hughes \cite{hughes} and O. Pironneau\cite{pironneau}. Later, in the nineties Finite Volume methods and Discontinuous Galerkin methods were proposed for non-solenoidal convective velocities (see for example A. Ern et al \cite{ern}.)

Recently we were faced with the problem of finding a good method for the computation of the probability density of a process via the Kolmogorov forward equation.  Here positivity and conservativity are essential. A more subjective criteria is the numerical diffusivity.  It became an opportunity to review the state of the art forty years after the above mentioned methods were proposed, what R. Glowinski would call a rear-guard battle. Nevertheless, the following methods are popular:
\begin{itemize}
\item The Primal Characteristic-Galerkin method (PCGM) proposed in  \cite{pironneau} is very precise but known to diverge in some cases when the viscosity is zero \cite{suli} and it is not conservative. It is convergent when mass-lumping is used \cite{optabata} but then it is too diffusive.
\item The Dual Characteristic-Galerkin method (DCGM) proposed in \cite{benque}  by J.P. Benque et al.  was never shown to converge except possibly when the initial and convected triangulations are intersected.
\item T.J.R. Hughes' streamline upwinding method (SUPG) \cite{hughes}, also called Galerkin Least-square upwinding \cite{johnson}, easy to implement, conservative and convergent but numerically  diffusive,  even when the \emph{upwinding parameter} is tuned to the problem.
\end{itemize} 

In the present note we study the DCGM with numerical quadrature for the nonlinear integral, prove that it is conservative, $L^2$-stable and convergent when the diffusion coefficient $\nu$ is not zero.  Proposition \ref{prop3}, below, shows that the method is $O(h+h^2/\delta t)$ when $\nu \gg h^2/\delta t$; $\delta t$ is the time step and $h$ is the  size of the edges of the triangulation.

The numerical section shows the superiority of DCGM over all 4 above cited methods. But DCGM is difficult to program. Indeed it is hard (but not computer intensive) to find in which element of the triangulation lies a given point, a well known problem of computational geometry \cite{shamos}.

Note also that the paper analyzes only the case of homogeneous Neumann condition. It ends with a numerical test with non-homogenous Dirichlet conditions for the Navier-Stokes equations, but the error analysis does not apply and it seems that it is numerically sensitive to the choice of the time step. 

\section{The Dual Characteristic-Galerkin Method}

Given a real parameter $\nu>0$,  a bounded open set $\Omega$ of $\R^d,~d=2,3$, a smooth velocity field $\va:~\Omega\times(0,T)\to \R^d$ and an initial condition $u^0: \R^d\to \R$, we wish to find $u:\Omega\times(0,T)\to \R$ such that, at all time $t\in(0,T)$, 
\eq{&&\label{contin}
\partial_t u + \va\cdot\nabla  u -\nu\Delta u =0,\quad u(0)=u^0\hbox{ in } \Omega,\quad \partial_n u=0\hbox{ on } \partial\Omega .%\hbox{ or } u=u^0  \hbox{ on }\Gamma.
}
Let $\bar\va$ be the extension of $\va$ by zero outside $\Omega$.  Define: $\dot{\bm\eta}(t)=\bar\va(\bm\eta(t)),~ \bm\eta(0)=\vx$ and $\bm\eta^\pm(\vx)=\bm\eta(\pm\delta t)$. Recall that 
\[
\partial_t u(\vx,t) + \va(\vx)\cdot\nabla  u(\vx,t) = \lim_{\delta t\to 0}\frac{1}{\delta t}[u(\vx,t)-u(\bm\eta^-(\vx),t-\delta t)].
\]

We assume that $\nabla\cdot\va=0$ and $\va\cdot{\bf n}=0$ at the boundary $\Gamma:=\partial\Omega$, so that $\bm\eta^\pm(\Omega)=\Omega$ and det$\nabla\bm\eta^\pm=1$. Hence two variational formulations of the problem discretized in time are feasible,
\eq{\label{varform}&&
\int_\Omega(\frac1{\delta t}(u^n\,\hat u - u^{n-1}\circ\bm\eta^-\,\hat u) + \nu\nabla u^n\cdot\nabla\hat u  )=0\quad \forall\hat u\in H^1(\Omega),\quad\add{(Primal~ form)},
\cr&&
\int_\Omega(\frac1{\delta t}(u^n\,\hat u - u^{n-1}\,\hat u\circ\bm\eta^+) + \nu\nabla u^n\cdot\nabla\hat u  )=0\quad \forall\hat u\in H^1(\Omega), \quad\add{(Dual~ form)}.
}
We have used $\bm\eta^+(\bm\eta^-(\vx))=\vx$  and,
\eq{\label{changevar}
\int_\Omega f(\vx)g(\bm\eta^-(\vx)) = \int_{\bm\eta^-(\Omega)} g(\vy)f(\bm\eta^+(\vy))/{\rm det}\nabla\bm\eta^-(\vy)=\int_\Omega g(\vy)f(\bm\eta^+(\vy)).
}
A spatial discretization with the Finite Element Method (FEM) of the first line in \eqref{varform} leads to the Primal Characteristic-Galerkin method (PCGM); on the second line it leads to the Dual Characteristic-Galerkin method (DCGM): finds $u^n\in V_{h}$ such that
%\color{blue}
\eq{\label{convectd}
\boxed{\int_\Omega \big(u^n_h \hat u_h + \delta t\nu\nabla u^n_h\cdot \nabla\hat u_h  \big)
= 
\sum_{i\in I} u^{n-1}_h(\bm\xi^i)\hat u_h(\bm\eta^i)\omega^i, \quad \forall \hat u_h\in V_{h},
}}
%\color{black} 
where,
\begin{itemize}
\item $\Omega$ is polygonal so as to be covered by a triangulation $\cup_k T^k$.
\item The points $\{\bm\xi^i\}_{i\in I}$ and \add{positive } weights 
$\{\omega^i\}_{i\in I}$ define a quadrature rule which must be exact  at least for  continuous piecewise-$P^2$ functions on the triangulation.  We assume that the quadrature is defined on triangles so as to write
\begin{equation}\label{quadra}
\sum_{i\in I}f(\bm\xi^i)\omega^i := \sum_k \sum_{i\in I(T^k)}f(\bm\xi^i)\omega^i_k,\quad I=\cup_k I(T^k).
\end{equation}
\begin{example}\label{ex1}
In 2D one may choose the quadrature points at the mid edges and $\omega^i_k=\frac13$, but more precise formulae are permitted.
\end{example}
\item $\bm\eta^i\in\Omega$ is an approximation of $\bm\eta^+$ with $|\bm\eta^i -\bm\eta^+(\bm\xi^i)|\le C\delta t^2.$ For example
%\eq{\label{eta}
%
\eq{\label{secondo}
\add{\bm\eta_a^+(\vx)} = \vx + \va(\vx)\delta t + \frac\sigma2\delta t^2\va(\vx)\cdot\nabla \va(\vx),\quad \sigma = 0 \hbox{ or } 1,\quad \add{\bm\eta^i=\bm\eta_a^+(\bm\xi^i)}.
}
\item $V_h$ is the $P^1$ continuous finite element space. 
\end{itemize}

\begin{proposition}
DCGM conserves mass in the sense that
\add{\[
\int_\Omega u^n_h =\int_\Omega u^0_h,~~\forall n.
\]}
\end{proposition}
\emph{Proof:}
Simply replace $\hat u_h$ by $1$ in the scheme.

\begin{proposition}\label{prop2}
\add{Assume that the triangulation is regular, in the sense of \cite{ciarlet}(p131) , i.e. for all triangles, the ratio of largest edge to the radius of the inscribed circle is bounded independently of $h$.} Then  DCGM is stable:
\[
\|u^n_h\|_{\nu\delta t} \le \left(1+|{\rm det}\underline{ {\bf A}}|\delta t^2+
C \frac{\add{h^2}}\nu \right)\|u_h^{n-1}\|_{\nu\delta t}
\]
where $\|v\|_{\nu\delta t}:=(| v|^2_0 + \delta t\nu|\nabla v|^2_0)^\frac12$, $C$ is a generic constant and $h$ is the  length of the longest edges in the triangulation.
\end{proposition}

\emph{Proof: }
The proof is given in 2D with the quadrature at the mid-edges (Example \ref{ex1}) and scheme \eqref{secondo}.

The discrete Cauchy-Schwarz inequality applied to the right hand-side of \eqref{convectd} combined with the choice $\hat u_h=u^n_h$ in \eqref{convectd}, leads to
\begin{eqnarray}\label{toto}&
\|u^n_h\|_{\nu\delta t}^2 &\le \left(\sum_{i\in I} u^{n-1}_h(\bm\xi^i)^2\omega^i\right)^\frac12\left(\sum_{i\in I}u^n_h(\bm\eta^i)^2\omega^i\right)^\frac12
%\cr&&
\le
\|u^{n-1}_h\|_{\nu\delta t}\left(\sum_{i\in I}u^n_h(\bm\eta^i)^2\omega^i\right)^\frac12,
\end{eqnarray}
because the quadrature is exact for $(u_h^{n-1})^2$ and because $|u_h^{n-1}|_0\le \|u^{n-1}_h\|_{\nu\delta t}$ . The map $\bm\xi\to\bm\eta_a^+(\bm\xi)$ defined by \eqref{secondo} transforms a triangle $T^k$ of the triangulation   into $\hat T^k$ and $\{\bm\eta^i,\omega^i\}_{i\in I}$ is a quadrature rule which is almost exact on $P^2$ functions of $\hat T^k$. 
We will show that, for some $C$,
\eq{\label{ineq}
\sum_k\sum_{i\in I(T^k)}u^n_h(\bm\eta^i)^2\omega^i_k
\le \left(1+ C (\add{\frac{h^2}\nu} + \delta t^2)\right)\|u^n_h\|^2_{\nu\delta t}. 
}

\subsubsection*{Proof of \eqref{ineq} in the linear case}
%%%%%%%%%%%%%%%%%%%%%
Assume that $\va$ is linear in ${\bf x}=(x,y)^T$  with $\nabla\cdot \va =0$, and consider \add{the case $\sigma=0$ in \eqref{secondo}},
\[
\bm\eta^+(\vx)={\bf x} +\delta t \va(\vx) = {\bf x} + \delta t\left[\begin{matrix} a^0_1 \cr a^0_2 \end{matrix}\right]
+
\delta t\left[\begin{matrix} \partial_x\va_1 x + \partial_y\va_1 y \cr \partial_x\va_2 x -\partial_x\va_1 y \end{matrix}\right] \add{= \bm\eta_a^+(\vx)}.
\]
It is not quite an isometry because det$ \nabla ({\bf x} + \va\delta t)= 1 -[(\partial_x\va_1)^2+\partial_y\va_1\partial_x\va_2]\delta t^2.$

Consider the quadrature at the mid edges with weight $\omega^i_k=\frac13|T^k|$, the area of $T^k$.
A triangle $(\vq^1,\vq^2,\vq^3)$ is transformed by $\bm\eta^+$ into the triangle $(\hat\vq^1,\hat\vq^2,\hat\vq^3)$ with 
\[
\hat\vq^j = \vq^j + \delta t\va^0 + \delta t(\nabla\va)^T\,\vq^j.
\]
Obviously a mid edge $\frac12(\vq^{j_1}+\vq^{j_2})$ of $T^k$ is mapped into a mid edge of $\hat T^k$.
Therefore, the only error is due to the variation of the area of the triangle: $|\hat T^k|=$det$\nabla(\vx+\delta t\va)|T^k|$. Indeed, as $u^n_h(\bm\eta^+)$ is affine on $T^k$ and because of \eqref{changevar},
\[
\sum_{i\in I(T^k)}u^n_h(\bm\eta^i)^2\omega^i_k =|(u_h\circ\eta^+)^2|_{0,\hat T^k}= (1-\delta t^2{\rm det}\nabla\va)|u_h^2|_{0,T^k},
\]
\add{because the quadrature is exact for $P^2$ functions;  $|f|_{0,T}$  is the integral of $f$ on $T$.}
\subsubsection*{Proof in the general case}
%%%%%%%%%%%%%%%%%%%%
Consider a triangle $T^k$ and a Taylor expansion of $\va$ about $\vx^0$, the center of $T^k$,
\begin{eqnarray*}&&
\va(\vx)=\va_0+\underline{ {\bf A}}(\vx-\vx^0) + \frac12(\vx-\vx^0)\otimes(\vx-\vx^0):\underline{\bm\Phi}(\vx).
\end{eqnarray*}
\add{With scheme \eqref{secondo}, for some bounded function $\bm\Psi$,}
\begin{eqnarray*}&&
\bm\eta_a^+(\vx)
= \vx + \delta t\Big(\va_0+\underline{\bf A}(\vx-\vx^0) + \frac12(\vx-\vx^0)\otimes(\vx-\vx^0):\underline{\bm\Phi}(\vx) 
\cr&&
\hskip 2.5cm  + \frac\sigma 2\delta t\big( 
\underline{\bf A}\va_0 
+ \underline{\bf A}^2(\vx-\vx^0)
\add{+  \va_0\otimes(\vx-\vx^0):\underline{\bm\Phi}(\vx)   
+ (\vx-\vx^0)\otimes(\vx-\vx^0) \bm\Psi(\underline{ {\bf A}},\underline{\bm\Phi},\nabla\underline{\bm\Phi}\big)\Big)}
\end{eqnarray*}
It is of the form
\[
\bm\eta_a^+(\vx)= \bm\eta_l(\vx) + \delta t(\vx-\vx^0)\otimes(\vx-\vx^0):\bm\Psi_1
\quad \hbox{ where }\quad \bm\eta_l(\vx) := \vx + \delta t( \va_1 + \underline{\bf A}_1(\vx-\vx^0)),
\]
and where $\va_1, \underline{\bf A}_1, \bm\Psi_1$ are affine in $\delta t$.

Recall the notation $\bm\eta^i:=\bm\eta_a^+(\bm\xi^i)$ and let $\bm\eta_l^i:=\bm\eta_l(\bm\xi^i)$. The segment $[\bm\eta_l^i,\bm\eta^i]$  cuts a finite number of edges of the triangulation. Let these intersections be $\{\bm\xi^i_j\}_1^{J-1}$. With the convention that $\bm\xi^i_0:=\bm\eta_l^i$ and $\bm\xi^i_J:=\bm\eta^i$, we can write 
\color{blue}
\[
u^n_h(\bm\eta^i)^2 - u^n_h(\bm\eta_l^i)^2= \sum_{0\le j\le J-1}(u^n_h(\bm\xi^i_{j+1})^2-u^n_h(\bm\xi^i_j)^2).
\]
Each term is continuously differentiable, so the following Taylor expansion is valid,
\[
u^n_h(\bm\eta^i)^2 - u^n_h(\bm\eta_l^i)^2= 2\sum_{0\le j\le J-1}u^n_h(\vx^i_j)\cdot\nabla u^n_j(\vx^i_j)(\bm\xi^i_{j+1}-\bm\xi^i_j)
\le
2\max_j|u^n_h(\vx^i_j)\cdot\nabla u^n_j(\vx^i_j)|\;|\bm\eta^i-\bm\eta_l^i|,
\]
where $\vx^i_{j}\in[\bm\xi^i_j,\bm\xi^i_{j+1}]$.
Let $\vx^i_M=\hbox{arg}\max_j|u^n_h(\vx^i_j)\cdot\nabla u^n_j(\vx^i_j)|$. Then we have found $\vx_M^i\in[\bm\eta^i,\bm\eta_l^i]$ such that,
\[
u^n_h(\bm\eta^i)^2 
\le 
u^n_h(\bm\eta_l^i)^2 
+ 2|u^n_h(\vx^i_M)\cdot\nabla u^n_j(\vx^i_M)|\;|\bm\eta^i-\bm\eta_l^i|.
\]
\color{black}
By hypothesis $\nabla\cdot\va=0$, so  $\underline{\bf A}$ is as above . Hence, $\vx\to\eta_l(\vx)$ being affine (see \eqref{toto}), 
$\sum_{i\in I(T^k)}u^n_h(\bm\eta_l^i)^2\omega^i_k $ is bounded by $(1-{\rm det}\underline{\bf A}_1\delta t^2)|u_h^n|^2_{0,T^k}$. Now  $|\bm\eta^i-\bm\eta_l^i|=\delta t(\bm\xi^i-\vx^0)\otimes(\bm\xi^i-\vx^0):\Psi_1|$, so,
\begin{eqnarray*}&&
\sum_{i\in I(T^k)}u^n_h(\bm\eta^i))^2\omega^i_k 
\le
(1-{\rm det}\underline{\bf A}_1\delta t^2)|u_h^n|^2_{0,T^k}+ h^2\delta t \|\bm\Psi_1\|_\infty
\sum_{i\in I(T^k)} 2 |u^n_h(\vx^i_M)\cdot\nabla  u^n_h(\vx^i_M)|\omega^i_k 
\end{eqnarray*}
A discrete Cauchy-Schwarz inequality leads to,
\begin{eqnarray*}&&
2|u^n_h(\vx^i_M)||\nabla u^n_h(\vx^i_M)|
\le
u^n_h(\vx^i_M)^2+|\nabla u^n_h(\vx^i_M)|^2
\le
\frac1{\nu\delta t}\left(u^n_h(\vx^i_M)^2+\nu\delta t|\nabla u^n_h(\vx^i_M)|^2\right).
\end{eqnarray*}
\color{blue}
At the cost of a multiplicative constant we may replace $\vx^i_M$ by $\bm\xi^{j(i)}$, the nearest quadrature point in the  triangle of $\vx^i_M$ and obtain, 
\begin{eqnarray*}&&
\sum_k\sum_{i\in I(T^k)}2|u^n_h(\vx^i_M)\cdot\nabla  u^n_h(\vx^i_M)|\omega^i_k 
%\cr&&
\le
\frac{C}{\nu\delta t}\sum_k\sum_{i\in I(T^k)}\left(u^n_h(\bm\xi^{j(i)})^2+\nu\delta t|\nabla u^n_h(\bm\xi^{j(i)})|^2\right)\omega^i_k
\le \frac{C'}{\nu\delta t}\|u^n_h\|^2_{\nu\delta t}.
\end{eqnarray*}
The last inequality holds for a regular triangulation because each quadrature point occurs at most $N$ times, finite,  and the $\omega^i_k$ differs from $\omega^{j(i)}_k$ at most by the ratio $R$ of areas of  triangles:
\begin{eqnarray*}&\ds 
\sum_{k,i\in I(T^k)}\left(u^n_h(\bm\xi^{j(i)})^2+\nu\delta t|\nabla u^n_h(\bm\xi^{j(i)})|^2\right)\omega^i_k
&\le
\sum_{k,i\in I(T^k)}\max\frac{\omega^{i}_k}{\omega^{j(i)}_k}\left(u^n_h(\bm\xi^{j(i)})^2+\nu\delta t|\nabla u^n_h(\bm\xi^{j(i)})|^2\right)\omega^{j(i)}_k
\cr&&
\le R\; N
\sum_{k,i\in I(T^k)}\left(u^n_h(\bm\xi^i)^2+\nu\delta t|\nabla u^n_h(\bm\xi^i)|^2\right)\omega^i_k .\end{eqnarray*}
\color{black}
 In the end,
\begin{eqnarray*}&&
\sum_k\sum_{i\in I(T^k)}u^n_h({\bm\eta}^i)^2\omega^i_k 
\le 
\left(1+|{\rm det}\underline{ {\bf A}}|\delta t^2+
C \frac{\add{h^2}}\nu \right)\|u^n_h\|^2_{\nu\delta t}.
\end{eqnarray*}
This proves \eqref{ineq} and completes the proof of Proposition \ref{prop2}.

\subsection{Error Estimates}
%%%%%%%%%%%%%%%%%%%%%%%%%%%%
%Let $\Pi_h$ be the $P^k$-interpolation operator. 
Let $u^n_e\in H^1(\Omega)$ be the solution of the continuous problem \eqref{contin} discretized in time and with the same $\bm\eta_a^+$ as in the discrete case; then let $u^n_{eh}\in V_h$ be the  projection of $u_e^n$ in the sense that
\eq{\label{projec}&&
\int_\Omega( u_e^n \hat u+\nu\delta t\nabla u_e^n\nabla \hat u )=\int_\Omega u_e^{n-1} \cdot\hat u\circ \bm\eta_a^+,
\quad \forall \hat u\in H^1(\Omega),,
\cr&&
\int_\Omega (u_{eh}^n \hat u_h+\nu\delta t\nabla u_{eh}^n\nabla\hat u_h )=\int_\Omega( u_e^n \hat u_h+\nu\delta t\nabla u_e^n\nabla \hat u_h ) \quad \forall \hat u_h\in V_h.
}
\begin{lemma}
Let $\epsilon^n_h=u^n_h-u^n_{eh}$ defined by \eqref{projec}. \add{Then,}
\eq{\label{recur}&&
\|\epsilon_h^n\|^2_{\nu\delta t}
\le 
\left(
1+ C (\add{\frac{h^2}\nu} + \delta t^2) 
\right)\|\epsilon^{n-1}_h\|_{\nu\delta t}^2
%\cr&&
+
C h^2\|\epsilon_h^{n-1}\|_{\nu\delta t}.
}
\end{lemma}

\emph{Proof }

 Let $Q$ be the quadrature \eqref{quadra},
\[
Q_\Omega(v,w) := \sum_{i\in I} v(\bm\xi^i)w(\bm\xi^i)\omega^i =\sum_k Q_{T^k}(v,w),
\quad
Q_{T^k}(v,w) = \sum_{i\in I(T^k)} v(\bm\xi^i)w(\bm\xi^i)\omega^i_k.
\]
 Then $\forall \hat u_h\in V_{h}$,
\eqn{&
\ds\int_\Omega \big(\epsilon^n_h \hat u_h + \delta t\nu\nabla \epsilon^n_h\cdot \nabla\hat u_h \big)
&= 
Q_\Omega( u_h^{n-1},\hat u_h \circ \bm\eta_a^+)
-\int_\Omega u^{n-1}_e\cdot \hat u_h \circ \bm\eta_a^+
\cr&&
= Q_\Omega( \epsilon_h^{n-1},\hat u_h \circ \bm\eta_a^+)
+ Q_\Omega( u_{eh}^{n-1},\hat u_h \circ \bm\eta_a^+)-\int_\Omega u^{n-1}_e\cdot \hat u_h \circ \bm\eta_a^+
}
Consequently
\eqn{&
\|\epsilon_h^n\|^2_{\nu\delta t}
&=Q_\Omega( \epsilon_h^{n-1},\epsilon_h^{n-1} \circ \bm\eta_a^+)
\cr&&
+ Q_\Omega( u_{eh}^{n-1}-u^{n-1}_e,\epsilon_h^{n-1} \circ \bm\eta_a^+)
\cr&&
+ Q_\Omega(u^{n-1}_e,\epsilon_h^{n-1} \circ \bm\eta_a^+) -\int_\Omega u^{n-1}_e\cdot \epsilon_h^{n-1}\circ \bm\eta_a^+.
}
A discrete Schwartz inequality is applied to the  first term on the right and then \eqref{ineq},
\[
Q_\Omega( \epsilon_h^{n-1},\epsilon_h^{n-1} \circ \bm\eta_a^+)
\le \left(
1+ C (\add{\frac{h^2}\nu} + \delta t^2) 
\right)\|\epsilon^{n-1}_h\|_{\nu\delta t}^2
\]
The second term is handled in the same way,
\eqn{&
Q_\Omega( u_{eh}^{n-1}-u^{n-1}_e,\epsilon_h^{n-1} \circ \bm\eta_a^+)
&\le 
\left(
1+ C (\add{\frac{h^2}\nu} + \delta t^2) 
\right)\|\epsilon_h^{n-1}\|_{\nu\delta t}\cdot \|u_{eh}^{n-1}-u^{n-1}_e\|_{0}
\cr&&
\le C h^2\left(
1+ C (\add{\frac{h^2}\nu} + \delta t^2) 
\right)\|\epsilon_h^{n-1}\|_{\nu\delta t}.
}
Finally the third term is bounded by the quadrature error on $\hat T^k$ for $u_e^{n-1}\circ\bm(\eta^+)^{-1}$,
\[
Q_\Omega(u^{n-1}_e,\epsilon_h^{n-1} \circ \bm\eta_a^+) -\int_\Omega u^{n-1}_e\cdot \epsilon_h^{n-1}\circ \bm\eta_a^+
\le
(1+C\delta t^2)h^2\|u_e^{n-1}\circ(\bm\eta_a^+)^{-1}\|_{3}\cdot\|\epsilon_h^{n-1}\|_{\nu\delta t}.
\]
Let us gather the pieces
\eq{&&
\|\epsilon_h^n\|^2_{\nu\delta t}
\le 
\left(
1+ C (\add{\frac{h^2}\nu} + \delta t^2) 
\right)\|\epsilon^{n-1}_h\|_{\nu\delta t}^2
%\cr&&
+
C h^2\|\epsilon_h^{n-1}\|_{\nu\delta t}
}

%%%%%%%%%
\begin{proposition}\label{prop3}
%%%%%%%%%%%%%%%%%%%%%%%%%%%%%%
\eq{\label{errorep}
\|\epsilon_h^n\|_{\nu\delta t}
\le
 \left(\|\epsilon^0_h\|_{\nu\delta t} + C \frac{h^2}{\delta t}\right)\left(1+C (\frac{h^2}{\nu }+ \delta t^2)\right)^n .
}
\end{proposition}
\emph{Proof}

Recurrence \eqref{recur} is of the type
\[
(\varepsilon^n)^2 - (\varepsilon^{n-1})^2 \le \alpha (\varepsilon^n)^2 + \beta\varepsilon^n
\]
with $\varepsilon^n=\|\epsilon_h^n\|_{\nu\delta t}$, $\beta=C h^2$ and $\alpha=C (\add{\frac{h^2}\nu} + \delta t^2)$.
It is rewritten as
\eqn{&&
\varepsilon^n - \varepsilon^{n-1} \le  \frac{\varepsilon^{n-1}}{\varepsilon^n +\varepsilon^{n-1}}(\alpha\varepsilon^{n-1} + \beta)\le  \alpha\varepsilon^{n-1} + \beta
\cr&&
\Rightarrow\quad
\varepsilon^n\le \varepsilon_0(1+\alpha)^n + C h^2\sum_{j=0}^{n-1}(1+\alpha)^j
\le
\varepsilon_0(1+\alpha)^n +\frac{(1+\alpha)^{n}-1}\alpha C h^2. 
}
The result derives from the fact that $n\le T/\delta t$ and  $(1+\alpha)^{n}-1\le n\alpha (1+\alpha)^{n-1}$.

\begin{remark}
Notice that the sequence is closed to the solution of the ODE in time $\varepsilon'=\frac1{2\delta t}(\alpha\varepsilon + \beta)$,
\[
\varepsilon(t) + \frac\beta\alpha  = (\varepsilon(0)+ \frac\beta\alpha)\exp(t\frac\alpha{2\delta t}),\hbox{ approximated by }\varepsilon(t)\approx\varepsilon(0)(1+t\frac\alpha{2\delta t}) + t\frac\beta{2\delta t} \hbox{ when }h^2<<\nu\delta t,
\]
\add{ because then $\frac{\alpha}{\delta t}<<1$}. So, at best, a tighter argument will only improve the constants in \eqref{errorep}.
\end{remark}
\begin{remark}
To derive the total error from $\epsilon^n_h$ is standard.  The time discretization being first order it produces and extra $O(\delta t)$ term , so the total error is of order $\delta t + \frac{h^2}\nu$, provided $h^2<\nu\delta t$. \add{Notice that here too, as for Primal Characterisic-Galerkin methods, $\delta t$ should not be chosen too small.}
\end{remark}

%%%%%%%%%%%%%%%%%%%%%%%%%
\section{Numerical Tests}
%%%%%%%%%%%%%%%%%%%%%%%%%%
\subsection{The Rotating Gaussian Bell}
%%%%%%%%%%%%%%%%%%%%%%%%%%%%%%
A point $\vx ^0=(\vx^0_1,\vx^0_2)^T$ convected by $ \va(\vx) = (-\vx_2,\vx_1)^T$ is in fact rotated at time $t$ to $ \vx^0(t)=(\vx^0_1\cos t + \vx^0_2\sin t , -\vx^0_1\sin t + \vx^0_2\cos t)^T$.  Consider
\eq{\label{exact}
u_e(\vx,t) = \frac {e^{-\frac{ r|\vx - \vx^0(t)|^2}{1+4\nu r t}}}{1+4\nu r t}
}
It verifies \eqref{contin} and $\partial_n u_e\approx 0$ if $r$ is large and $\nu$ is small.

A Delaunay-Voronoi mesh generator is used for the triangulations of the unit circle.  We tested 3 meshes with 926, 3601 and 14071 vertices, corresponding respectively to $N=100$, 200 and 400 boundary vertices.  The corresponding number of time steps chosen are  33, 66 and 133.

The other parameters are  $\vx^0_1=0.35$, $\vx^0_2=0$,$ T=2\pi$, $\nu=10^{-4}$ or $0.01$, $r=10$.

\subsection{Convergence Study}
%%%%%%%%%%%%%%%%%%%%%%%%%%%%%
In this section $\nu=10^{-4}$.

The differential equation is discretized by \eqref{secondo} with $\sigma=1$.
$V_h$ is constructed with the linear continuous triangular finite element method and the nonlinear integral is approximated with the mid-edges as quadrature points of Example \ref{ex1} or a 9-points quadrature per triangle \cite{FF}. 

\begin{figure}[http]
%%%%%%%%%%%%%%%%%%%%
\begin{center}
\begin{minipage} [b]{0.45\textwidth}
\begin{tikzpicture}[scale=0.8]
\begin{axis}[legend style={at={(1,1)},anchor=north east}, compat=1.3,
   xlabel= {Number of vertices},
  ylabel= {$L^2$-error},
  xmode=log, ymode=log,
  ]
\addplot[thick,solid,color=blue,mark=*, mark size=1pt] table [x index=0, y index=1]{fig/dualbell2.txt};
\addlegendentry{9 pt quad}
\addplot[thick,solid,color=red,mark=none, mark size=1pt] table [x index=0, y index=2]{fig/dualbell2.txt};
\addlegendentry{3 pt quad}

\addplot[thick,dotted,color=black,mark=o, mark size=1pt] table [x index=0, y index=4]{fig/dualbell2.txt};
\addlegendentry{error$\sim h$}
\end{axis}
\end{tikzpicture}
\caption{\label{error1}  Plot (log-log scales) of $L^2$ error versus vertices number and effect of quadratures on the precision. }
\end{minipage}
\hskip 1cm
\begin{minipage} [b]{0.45\textwidth} % errorkoba.txt
\includegraphics[width=\textwidth]{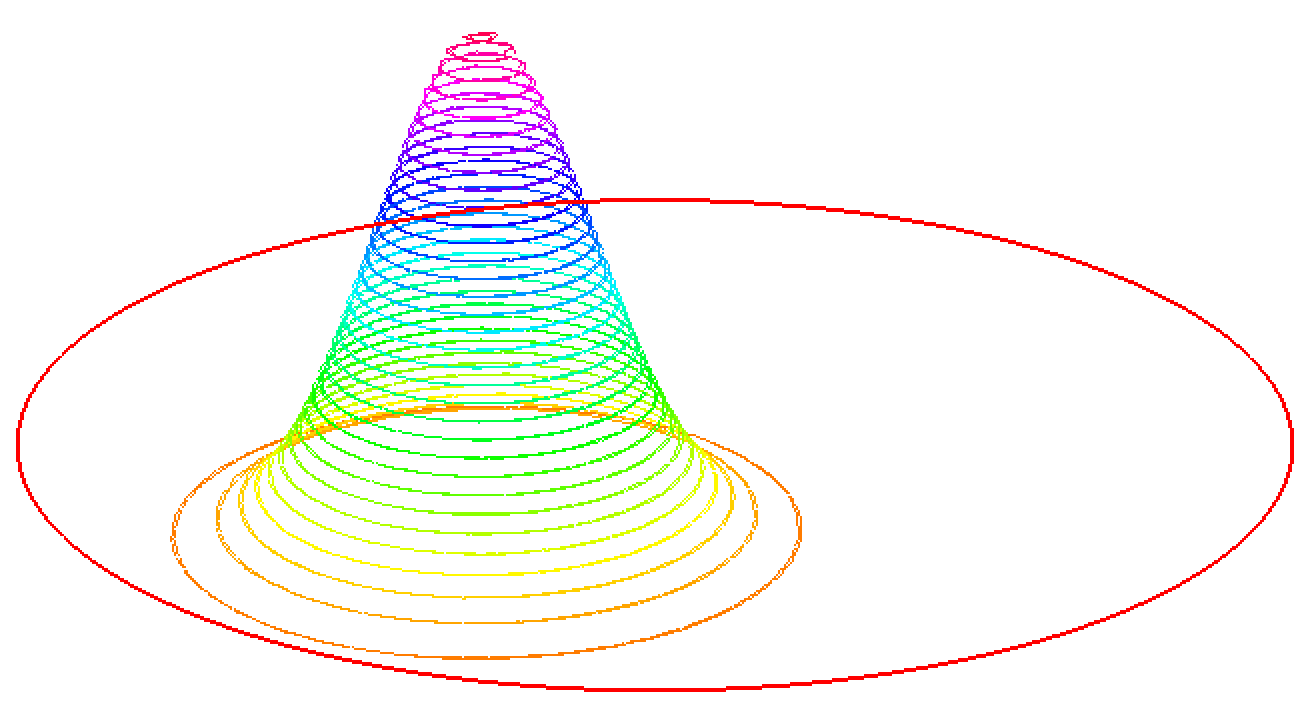}
\caption{\label{errorkoba} Gaussian Bell after one turn and exact solution.  The level lines of both surfaces are very near to each others. Level lines values are as in Fig. \ref{convectone}. }
\end{minipage}
\end{center}
\end{figure}

Figure \ref{error1} shows the convergence rate  and Figure \ref{errorkoba} shows the Gaussian bell after one turn. It is difficult to see the difference with the exact solution. 

\begin{figure}[htp]
%%%%%%%%%%%%%%%%%%%%
\begin{center}
\begin{minipage} [b]{0.5\textwidth}
\includegraphics[width=\textwidth]{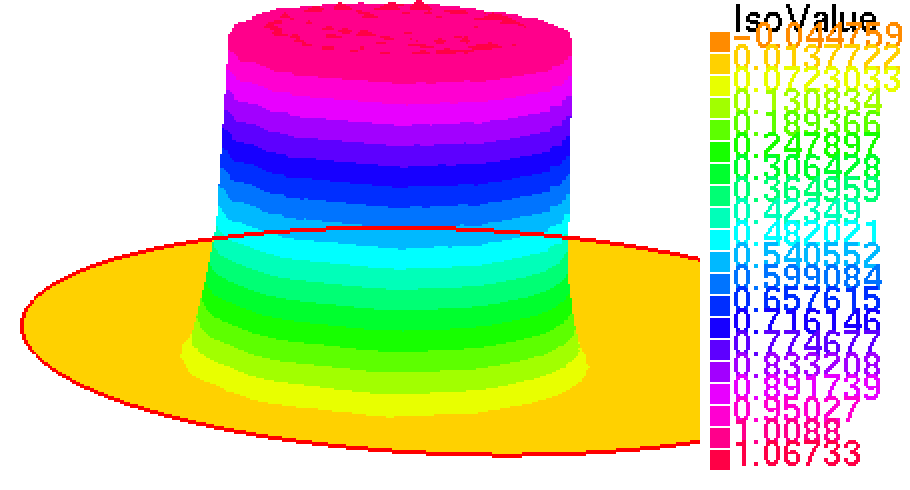
}
\caption{\label{convectone} $u^0={\bf 1}_{(x-0.3)^2+y^2<0.15}$and $u_h^T$ after one turn. Notice there is almost no oscillation and no numerical diffusion. }
\end{minipage}
\hskip 1cm
\begin{minipage} [b]{0.35\textwidth} % errorkoba.txt
\includegraphics[width=\textwidth]{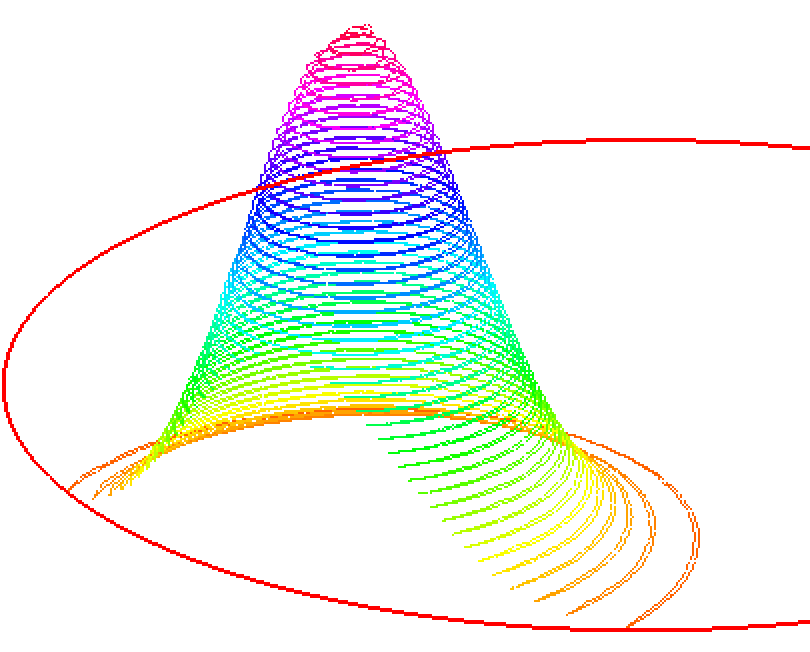}
\caption{\label{errorkobb} Gaussian bell crossing the boundary, because initially $x_0=0.5$, after one turn and exact solution.}
\end{minipage}
\end{center}
\end{figure}

A discontinuous function is subject to the rotating field to test the robustness with respect to discontinuity.  Results are on Figure \ref{convectone}.
Finally, as shown by Figure \ref{errorkobb} $u_h$ need not be zero at the boundary. 
Figures \ref{errorkoba}, \ref{convectone} and \ref{errorkobb}  have been computed with $N=200$.
Table \ref{pcc} shows the positivity and conservativity of the method.
\begin{table}[htp]
\caption{Positivity, Conservativity and Convergence}
\begin{center}
\begin{tabular}{|c||c|c|c|c|}
$N$ & $\min u_h$ & $\max u_h$ &$\int_\Omega u_h$ & $L^2$-error
\cr\hline\cr
100& -1.13689e-08& 0.643741 &0.156945 &0.0112869\cr
200& 1.94281e-11& 0.664612& 0.156998& 0.00282539\cr
400& 1.94281e-11& 0.665645& 0.156962& 0.000763338\cr
Exact &1.94281e-11& 0.665268& 0.156965& 0\cr
\end{tabular}
\end{center}
\label{pcc}
\end{table}%

\section{Comparison with other methods}
In this section $\nu=0.01$ and by default $N=200$.

We ran the same tests with 4 other popular methods: PCGM \cite{pironneau}, SUPG \cite{hughes}, DDG \cite{ern} and no upwinding \cite{arak}.
Streamline Upwinding Galerkin (SUPG) reads:
\[
\int_\Omega(\frac{u^n_h-u^{n-1}_h}{\delta t} + \va\cdot\nabla u )(w_h+\alpha \va\cdot\nabla w_h) 
				+\int_\Omega\nu \nabla u^n_h\cdot\nabla w_h= 0
\]
for all $w_h\in V_{h}$;  $\alpha=0.3$ in the numerical test.

With homogeneous Dirichlet conditions the Dual Discontinuous-Galerkin (DDG) methods is:
\[
\int_\Omega((\frac{u^n_h-u^{n-1}_h}{\delta t}+\va\cdot\nabla u^n_h)w_h +\nu\nabla u^n_h\cdot\nabla w_h)
  +\int_E w_h(\alpha | \vn\cdot \va |- \frac12 \vn\cdot \va) [u^n_h]=0
\]
for all $w_h\in V_{h}$; $\alpha=0.5$ in the numerical test.
Here $E$ is the set of inner edges and $[b]$ is the jump of $b$ across an edge of $E$.

Finally the centered method which keeps the convective terms as is
\[
\int_\Omega((\frac{u^n_h-u^{n-1}_h}{\delta t}+\va\cdot\nabla u^n_h)w_h +\nu\nabla u^n_h\cdot\nabla w_h)
  =0 \quad \forall w_h\in V_{h}\, .
\]
  A CFL condition $\delta t \le c(\nu) h^2$ is necessary for stability, so the method is not viable for small $\nu$. 
  
\begin{figure}[htp]
%%%%%%%%%%%%%%%%%%%%
\begin{center}
\begin{minipage} [b]{0.3\textwidth}
\begin{tikzpicture}[scale=0.6]
\hskip -0.5cm
\begin{axis}[legend style={at={(1,1)},anchor=north east}, compat=1.3,
   xmax=1.2,
   xlabel= {x},
  ylabel= {$u_h(x,0)$},
  ]
\addplot[thick,solid,color=blue, mark size=1pt] table [x index=0, y index=1]{fig/cut100.txt};
\addlegendentry{PCGM}
\addplot[thick,solid,color=red, mark size=1pt] table [x index=0, y index=2]{fig/cut100.txt};
\addlegendentry{DCGM}
\addplot[thick,solid,color=green, mark size=1pt] table [x index=0, y index=3]{fig/cut100.txt};
\addlegendentry{SUPG}
\addplot[thick,solid,color=magenta, mark size=1pt] table [x index=0, y index=4]{fig/cut100.txt};
\addlegendentry{DDG}
\addplot[thick,dashed,color=black, mark size=1pt] table [x index=0, y index=5]{fig/cut100.txt};
\addlegendentry{Centered}
\addplot[thick,dotted,color=black, mark=+,mark size=1pt] table [x index=0, y index=6]{fig/cut100.txt};
\addlegendentry{Exact}
\end{axis}
\end{tikzpicture}
\end{minipage}
\hskip 0.5cm
\begin{minipage} [b]{0.3\textwidth} % errorkoba.txt
\begin{tikzpicture}[scale=0.6]
\begin{axis}[legend style={at={(1,1)},anchor=north east}, compat=1.3,
   xmax=1.2,
   xlabel= {x},
  ]
\addplot[thick,solid,color=blue, mark size=1pt] table [x index=0, y index=1]{fig/cut200.txt};
\addlegendentry{PCGM}
\addplot[thick,solid,color=red, mark size=1pt] table [x index=0, y index=2]{fig/cut200.txt};
\addlegendentry{DCGM}
\addplot[thick,solid,color=green, mark size=1pt] table [x index=0, y index=3]{fig/cut200.txt};
\addlegendentry{SUPG}
\addplot[thick,solid,color=magenta, mark size=1pt] table [x index=0, y index=4]{fig/cut200.txt};
\addlegendentry{DDG}
\addplot[thick,dashed,color=black, mark size=1pt] table [x index=0, y index=5]{fig/cut200.txt};
\addlegendentry{Centered}
\addplot[thick,dotted,color=black, mark=+,mark size=1pt] table [x index=0, y index=6]{fig/cut200.txt};
\addlegendentry{Exact}
\end{axis}
\end{tikzpicture}
\end{minipage}
\hskip 0.5cm
\begin{minipage} [b]{0.3\textwidth} % errorkoba.txt
\begin{tikzpicture}[scale=0.6]
\begin{axis}[legend style={at={(1,1)},anchor=north east}, compat=1.3,
   xmax=1.2,
   xlabel= {x},
  ]
\addplot[thick,solid,color=blue, mark size=1pt] table [x index=0, y index=1]{fig/cut400.txt};
\addlegendentry{PCGM}
\addplot[thick,solid,color=red, mark size=1pt] table [x index=0, y index=2]{fig/cut400.txt};
\addlegendentry{DCGM}
\addplot[thick,solid,color=green, mark size=1pt] table [x index=0, y index=3]{fig/cut400.txt};
\addlegendentry{SUPG}
\addplot[thick,solid,color=magenta, mark size=1pt] table [x index=0, y index=4]{fig/cut400.txt};
\addlegendentry{DDG}
\addplot[thick,dashed,color=black, mark size=1pt] table [x index=0, y index=5]{fig/cut400.txt};
\addlegendentry{Centered}
\addplot[thick,dotted,color=black, mark=+,mark size=1pt] table [x index=0, y index=6]{fig/cut400.txt};
\addlegendentry{Exact}
\end{axis}
\end{tikzpicture}
\end{minipage}
\caption{\label{compare}  Plot of $x\to u_h(x,0)$ computed by the 5 methods, at $N=100$ (left), $N=200$ (middle) and $N=400$ (right) . }
\end{center}
\end{figure}
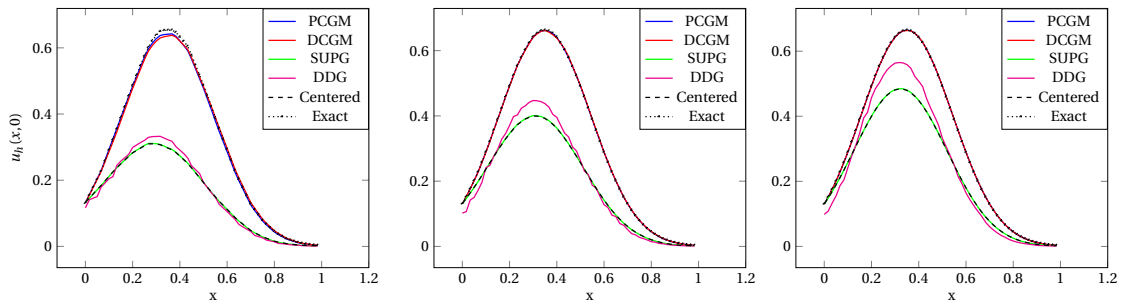

\begin{figure}[htp]
%%%%%%%%%%%%%%%%%%%%
\begin{center}
\begin{minipage} [b]{0.3\textwidth}
\includegraphics[width=1.\textwidth]{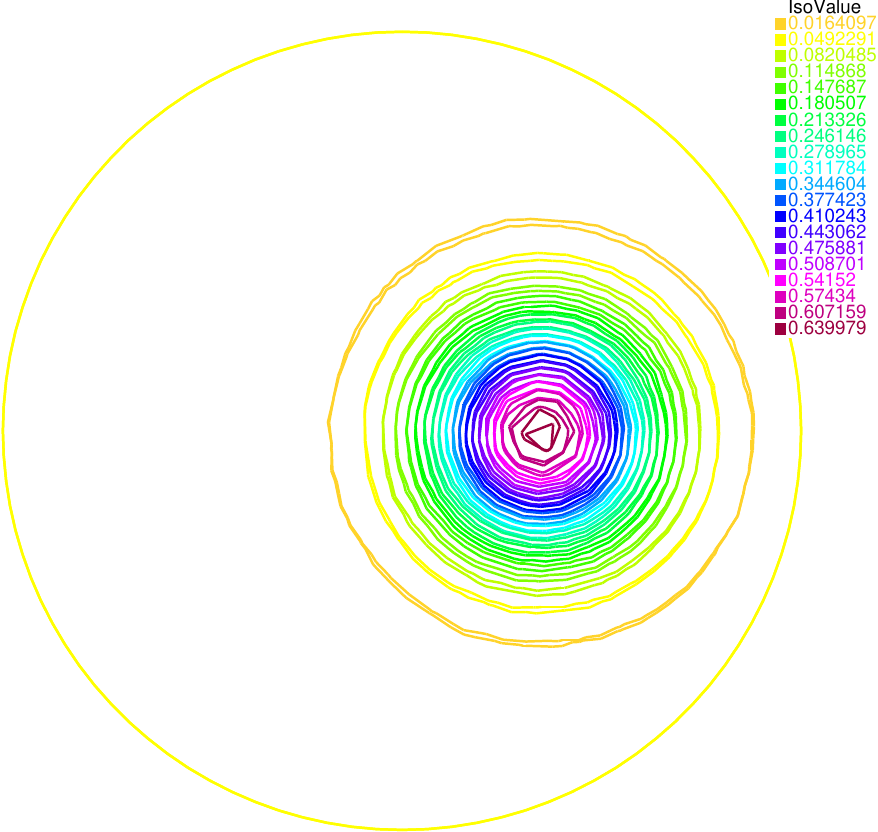}
\caption{\label{PCGM}  Bell computed with $N=100$ and  with PCGM after one turn and exact solution (level lines are essentially on top of each other). }
\end{minipage}
\hskip 0.5cm
\begin{minipage} [b]{0.3\textwidth} % errorkoba.txt
\includegraphics[width=1.\textwidth]{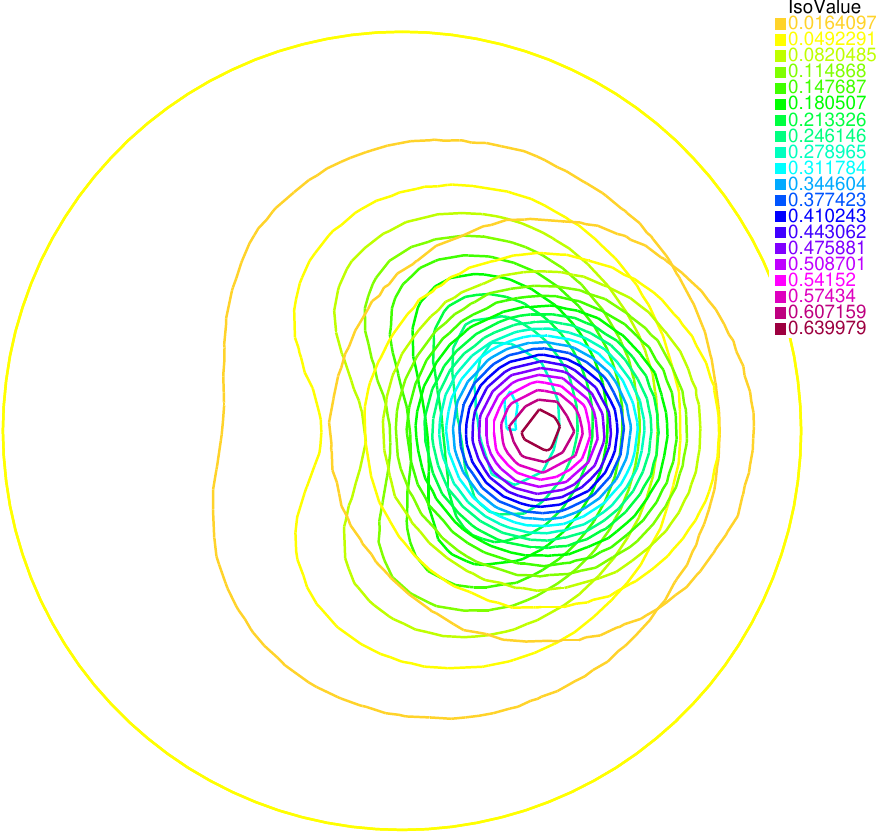}
\caption{\label{SUPG}  Bell computed with $N=100$ and  with SUPG after one turn and exact solution. Phase error, flatness error and maximum error are visible. }
\end{minipage}
\hskip 0.5cm
\begin{minipage} [b]{0.3\textwidth}
\includegraphics[width=1.\textwidth]{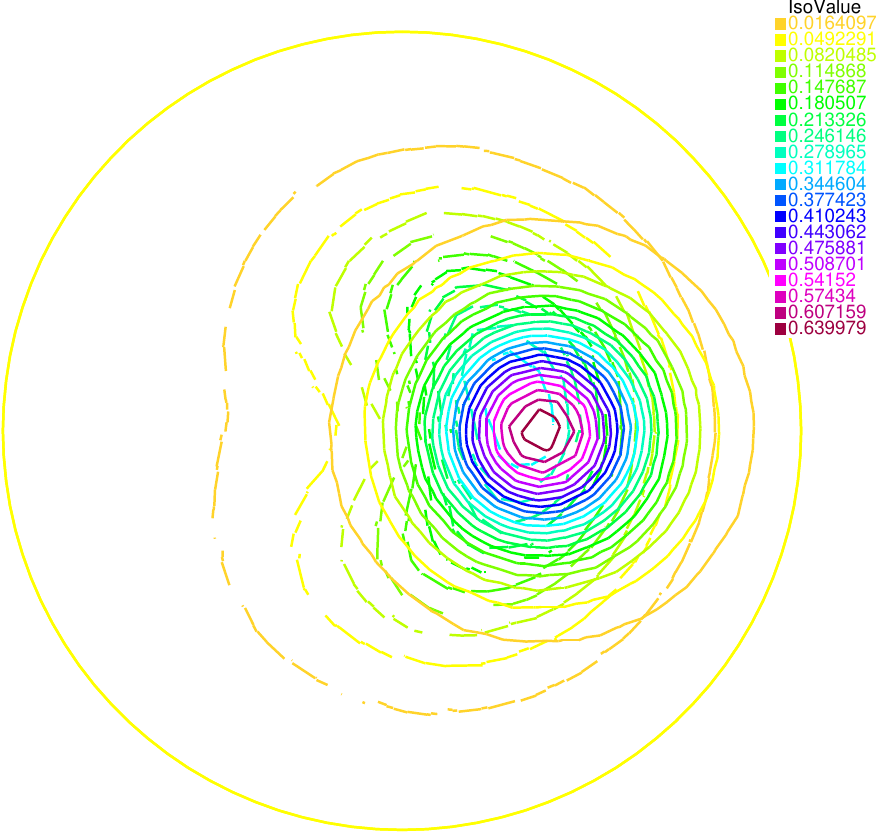}
\caption{\label{DDG} Bell computed with $N=100$ and  with DDG elements after one turn and exact solution. Phase, flatness and maximum error are visible.}
\end{minipage}
\end{center}
\end{figure}

\begin{figure}[htp]
%%%%%%%%%%%%%%%%%%%%
\begin{center}
\begin{minipage} [b]{0.45\textwidth}
\begin{tikzpicture}[scale=0.65]
\begin{axis}[legend style={at={(1,1)},anchor=north west}, compat=1.3,
   xlabel= {$N$, the number of boundary points},
  ylabel= {$L^2$-error},
  xmode=log, ymode=log,
  ]
\addplot[thick,solid,color=blue,mark=*, mark size=1pt] table [x index=0, y index=4]{fig/errorconvp.txt};
\addlegendentry{PCGM}
\addplot[thick,solid,color=red,mark=none, mark size=1pt] table [x index=0, y index=4]{fig/errorconvd.txt};
\addlegendentry{DCGM}
\addplot[thick,dotted,color=black,mark=+, mark size=1pt] table [x index=0, y index=4]{fig/errorsupg.txt};
\addlegendentry{SUPG}
\addplot[thick,dotted,color=red,mark=o, mark size=1pt] table [x index=0, y index=4]{fig/errordc.txt};
\addlegendentry{DC}
\addplot[thin,dashed,color=blue,mark=*, mark size=1pt] table [x index=0, y index=4]{fig/errorcenter.txt};
\addlegendentry{Centered}
\addplot[thin,dashed,color=black,mark=none, mark size=1pt] table [x index=0, y index=1]{fig/errh.txt};
\addlegendentry{error$\sim h$}
\end{axis}
\end{tikzpicture}
\caption{\label{converge2}  Plot (log-log scales) of $L^2$ error versus $N$. Both characteristic methods are equally precise and the other methods (SUPG, DDG, no upwding) are equally coarse. }
\end{minipage}
\hskip 0.5cm
\begin{minipage} [b]{0.45\textwidth}
\includegraphics[width=0.7\textwidth]{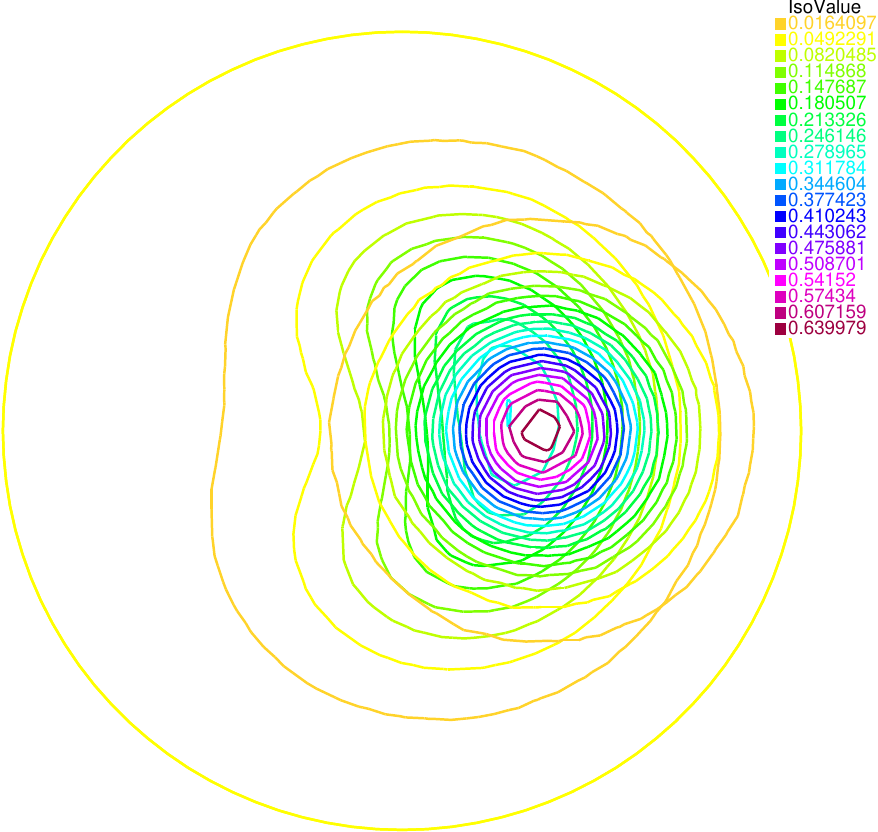}
\caption{\label{centered} Bell computed with $N=100$ and  with the centered FEM (i.e. without upwinding). There are ten times more time steps to perform a turn. Phase error, maximum error and flatness error are visible.}
\end{minipage}
\end{center}
\end{figure}

  Figure \ref{compare} shows the horizontal cross sections of the Gaussian bell in the $x$ direction after one turn for all 5 methods.  Obviously PCGM and DCGM perform better, with the advantage that DCGM is convervative and convergence is proved.  The level lines of the Gaussian bell after one turn are shown on Figures \ref{PCGM}, \ref{SUPG}, \ref{DDG} and \ref{centered} and the positivity and conservativity on Table \ref{tab2}.
Finally the convergence rates are shown in Figure \ref{converge2}.

%%%%%%%%%%%%%%%%%%%%%%%%
%\begin{figure}[htp]
%%%%%%%%%%%%%%%%%%%%%
%\begin{center}
%\begin{minipage} [b]{0.45\textwidth}
%\includegraphics[width=\textwidth]{fig/fullbell.png}
%\caption{\label{converge2} Gaussian bell computed with Primal Characteristics after one turn and exact solution (level lines are esstentially on top of each other). }
%\end{minipage}
%\hskip 1cm
%\begin{minipage} [b]{0.45\textwidth} % errorkoba.txt
%\includegraphics[width=\textwidth]{fig/bdyeffect.png}
%\caption{\label{errorkoba} Gaussian bell computed with Primal Characteristics after one turn and exact solution .}
%\end{minipage}
%\end{center}
%\end{figure}

%%%%%%%%%%%%%%%%%%%%%%%

\begin{table}[htp]
\caption{Comparison of the methods at N=200 after one turn.}
\begin{center}
\begin{tabular}{c||c|c|c|c}
Method & $\min u_h$ & $\max u_h$ &$\int_\Omega u_h$ & $L^2$-error\cr
\hline
$u_e$ intorpolated & 1.94281e-11& 0.66339& 0.156984 &\cr
PCGM & 1.94281e-11& 0.662813& 0.156777& 0.00277886\cr
DCGM & 1.94281e-11& 0.664612& 0.156998& 0.00282539\cr
SUPG & 1.94281e-11& 0.40193& 0.157103& 0.0893023\cr
DDG & 2.27941e-06& 0.448727& 0.157102& 0.0847009\cr
Centered & 1.94281e-11& 0.400491& 0.157099& 0.0894042\cr
\end{tabular}
\end{center}
\label{tab2}
\end{table}%

\section{Application to the Kolmogorov Equation for Heston's Model}
%%%%%%%%%%%%%%%%%%%%%%%%%%%%%%%%%%%%%%%%%%%%%%%%%%%%%%%%%%%%%%%%%%%
Let $\E[f]$ be the expected value of a random  $f$.  In quantitative finance Heston's model \cite{heston} is,
\eq{&&
{d X_t = X_t(rd t + \sqrt{Y_t}d W^1_t), 
\quad
d Y_t = \kappa(\theta-Y_t)d t + \lambda\sqrt{Y_t}d W^2_t},
\cr&&
\E[d W^1_td W^2_t]=\rho,\quad ~~X_0=\N(\mu,\sigma)
,\quad ~~Y_0=\N(\mu',\sigma').
}
It is popular to set the (undiscounted) price of a ``Put" to be $P_T=\E(K-X_T)^+$ at time $T$  where $K$ is the ``strike". Here   the random process $t\to \{X_t,Y_t\}$ is driven by its initial conditions $\{X_0,Y_0\}$ and the two normal Brownian motions $t\to W^i_t,~i=1,2$ with correlation $\rho$. The initial conditions are  Gaussian random variables of means $\mu,\mu'$ and standard deviations $\sigma,\sigma'$. The parameters $r,\kappa,\theta$ and $\lambda$ are positive real numbers.
Kolmogorov's theorem gives  the PDF $u\in L^2(\R_+^2)$ of  $\{X_t,Y_t\}$:  for all $\{x,y,t\}\in\R_+^2\times(0,T)$,
\eq{&&
{\small \partial_t u +\nabla\cdot\left[\begin{matrix}r x u\cr \kappa(\theta-y)u\end{matrix}\right] -\nabla^2:\left(\left[\begin{matrix} x^2y & \lambda x y\cr \lambda xy & \lambda^2 y\end{matrix}\right]\frac u2\right)=0,}
\qquad
 u_{|t=0} = G_{\mu,\sigma}(x)G_{\mu',\sigma'}(y),
}
where $G$ is the Gaussian curve.  Then {$ P_T=\int_{\R_+^2}(K-x)_+ u_T(x,y)$}.  Computing $P_T$ for large $T$ is a challenge because it is essential to keep having $\int_{\R_+^2}u_t=1$ for all $t$ and $u(x,y)\ge 0$ for all $x\ge0, y\ge 0$.

\begin{figure}[htbp]
\begin{center}
\includegraphics[width=\textwidth]{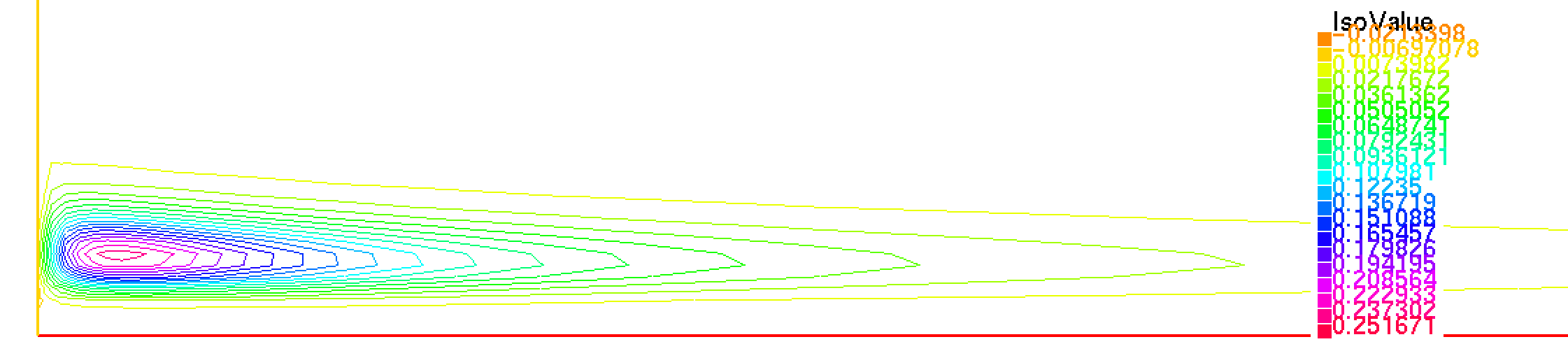}
\caption{The level lines of the PDF of Heston's model at time T=10.}
\label{hestonf}
\end{center}
\end{figure}

We computed $u_T$ at $T=10$ with DCGM when $r=0.03$, $K=75$, $\mu=50$, $\kappa=2$, $\theta=0.1$, $\lambda=0.2$, $\rho=-0.5$, $\mu'=0.75$, $\sigma=10$, $\sigma'=0.1$.  The results are in Figure \ref{hestonf} after 1500 time iterations and a mesh of $150\times 150$ vertices.
No negative values are observed and by construction $\int_{\R_+^2}u=1$. 

\section{Non Homogeneous Dirichlet Conditions}

\begin{figure}[htbp]
\begin{center}
\includegraphics[width=\textwidth]{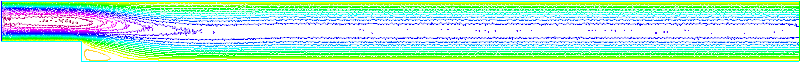}
\caption{Stationary solution of the Navier-Stokes equation at Reynold $50$. The level lines of the horizontal component of the fluid velocity are shown. The color scale is the same as that of Figure \ref{convectone}. The size of the recirculation is 3 times the height of the step as expected \cite{thomasset}.}
\label{nsres}
\end{center}
\end{figure}

Equation \eqref{changevar} is wrong when $\va\cdot\vn|_{\Gamma}\ne 0$. To compensate with the fact that $\eta^-(\Omega)\ne\Omega$,  a correction must be added (resp. subtracted)  
outside (resp. inside)  $\Gamma$ if $\va\cdot\vn|_{\Gamma}$ is negative (reps. positive). For Dirichlet conditions $u=u_\Gamma$, we propose to replace \eqref{convectd} by: find $u^n_h-u_\Gamma\in V_{0h}$ such that
\eq{\label{convectdd}
\int_\Omega \big(u^n_h \hat u_h + \delta t\nu\nabla u^n_h\cdot \nabla\hat u_h  \big)
-\int_{\Gamma}\delta t\va\cdot\vn u^n_h \hat u_h
= 
\sum_{i\in I} u^{n-1}_h(\bm\xi^i)\hat u_h(\bm\eta^i)\omega^i, \quad \forall \hat u_h\in V_{0h},
}
This formulation was tested on the Navier-Stokes equations for the backward step problem, using the $P^2-P^1$ element. Results are on Figure \ref{nsres}.  However the results are better without the boundary integral on right, \add{so something is afoot, the problem is open.}  

\bibliographystyle{plain}
\bibliography{references}

\end{document}